\documentclass[a4paper,12pt]{amsart}
\usepackage{times,a4wide,mathrsfs,amsthm,amssymb}
 \usepackage[cal=boondoxo]{mathalfa} 
\usepackage{hyperref}
\usepackage{color}
\usepackage[all]{xy}

\newtheorem{nonumberingt}{Acknowledgements}

\newtheorem{thm}{Theorem}[section]
\newtheorem*{thm*}{Theorem}
\newtheorem*{corr*}{Corollary}
\newtheorem{lemma}[thm]{Lemma}

\newtheorem{prop}[thm]{Proposition}
\newtheorem{corr}[thm]{Corollary}
\theoremstyle{definition}

\newtheorem*{notation}{Notation}

\theoremstyle{remark}
\newtheorem{conjecture}[thm]{Conjecture}

\newtheorem{rmk}[thm]{\textit{Remark}}
\renewcommand{\proof}{\noindent\textit{Proof}\/: \,\,}
\newcommand{\Q}{{\mathbf{Q}}}

\newcommand\LL{{\mathcal L}} 
\newcommand{\comp}{\raise1pt\hbox{{$\scriptscriptstyle\circ$}}}
 \def\im{\operatorname{Im}}
  \def\ker{\operatorname{Ker}}
\def\lset{\{}  
\def\rset{\}}  
 
\def\st{\mid}   
\def\sett#1#2{\lset #1 \st #2 \rset}

\def\End{\mathop{\rm End}\nolimits}
\def\Corr{\mathop{\rm Corr}\nolimits}

\def\chow{\mathop{\mathsf{CH}}\nolimits}

\def\id{\text{\rm id}}
\def\mapright#1{\mathop{\vbox{\ialign{
                                ##\crcr
    ${\scriptstyle\hfil\;\;#1\;\;\hfil}$\crcr
 \noalign{\kern2pt\nointerlineskip}
    \rightarrowfill\crcr}}\;}}

\def\mapleft#1{\mathop{\vbox{\ialign{
                                ##\crcr
    ${\scriptstyle\hfil\;\;#1\;\;\hfil}$\crcr
 \noalign{\kern2pt\nointerlineskip}
    \leftarrowfill\crcr}}\;}}

\newcommand\rarrow[3]{\smash{\mathop{\hbox to#3{\rightarrowfill}}\limits
^{\scriptstyle#1}_{\scriptstyle#2}}}
 
\newcommand\larrow[3]{\smash{\mathop{\hbox to#3{\leftarrowfill}}\limits
^{\scriptstyle#1}_{\scriptstyle#2}}}

\def\into{\hookrightarrow}
\def\onto{\twoheadrightarrow}
\begin{document}
\date{4 October    2017} 
  \author[Chris Peters]
{Chris Peters}
\address{Technical University Eindhoven, Netherlands.}
\email{c.a.m.peters@tue.nl}

\title{On a  motivic interpretation of  primitive,  variable and fixed cohomology}

 \maketitle

\section{Introduction}

This note aims to address the motivic nature of some classical cohomological results of Lefschetz. The first is the
 Lefschetz decomposition of the cohomology of a smooth projective manifold.    The second is a consequence of
Lefschetz'  hyperplane theorem, namely the 
splitting   of the cohomology of   a complete intersection  into a summand which comes from the surrounding variety, the "fixed part",
  and a supplementary summand, the "variable" part.  Explicitly, fix an $(d+r)$-dimensional projective manifold $M$ and an ample
  line bundle $\LL$ on $M$;  let  $X=H_1\cap \cdots\cap H_r$   be a smooth  complete intersection of $r$ divisors $H_j\in |\LL|$, $j=1,\dots,r$  and let 
   $i:X\into M$ be  the inclusion.
 With
  \begin{equation}\label{eqn:FixAndVarCoh}
  \aligned
  H^d(X)_{\rm fix}  &:=    \im( i^* : H^{d}(M) \to H^d(X)) \\
  H^d(X)_{\rm var} &:= \ker(i_*:H^d(X)\to H^{d+2r}(M))
  \endaligned
  \end{equation}
there is an orthogonal direct sum decomposition
 \begin{equation}\label{eqn:FixAndVarCoh2}
    H^d(X)=  H^d(X)_{\rm fix}\oplus  H^d(X)_{\rm var} .
\end{equation} 

In general it seems hard to show the motivic nature of these results and some conditions will be needed,
Clearly, a first ingredient one needs is the existence of a correspondence inducing the inverse of the Lefschetz operator on $H^*(M)$.
This is Lefschetz' conjecture $B(M)$. The second comes from a concept introduced by 
Kimura \cite{Kim} and O'Sullivan, the concept of finite-dimensionality
for motives.  They conjecture  that all motives are finite-dimensional.
 The main   result  of this note  is that  the primitive  decomposition  for the cohomology of $M$  as well as   the splitting \eqref{eqn:FixAndVarCoh2} is motivic 
 is  motivic provided  these two conjectures hold for $M$.  \footnote{ For the comfort   of the reader   some facts about Chow motives are placed together 
  in Section~\ref{sec:hlp}.} In fact, only a consequence of finite dimensionality is used, namely a certain nilpotency result which is stated as \eqref{eqn:Nilp}.
  
  It is known that both Kimura's conjecture and conjecture $B$ are verified  for example for  $M$ a  projective space, or an abelian variety.   
  For these examples   the motive of  $M$ is well understood and the primitive decomposition is probably well known. 
  See  e.g. Diaz'  explicit results   \cite{diaz} for abelian varieties.
  The motivic  nature of the splitting \eqref{eqn:FixAndVarCoh2}  for  complete intersections  $X\subset M$ 
  shows that   the relevant motivic information is hidden in the variable motive. 
 
  This can be taken advantage of in situations where the motive of $M$ is too large.    Let me illustrate this  with  the 
  Bloch conjecture \cite{B} for surfaces. Recall that the latter     states that if $p_g =0$,    for zero-cycles   homological equivalence and  rational  equivalence coincide  so that
  $\chow_0 $ is "small".
 In the present setting, assuming that one has a \emph{complex} complete intersection surface
 $X\subset M$, such that  $h^{2,0}(M)\not=0$, then, by Lefschetz' theorem on hyperplane sections $h^{2,0}(X)\not=0 $,  
 and then, by a result of Mumford \cite{M},
  the Chow group  of zero cycles on $X$  is huge. However, it may happen that the variable submotive of $X$, or a smaller submotive thereof
  does satisfy the conditions for Bloch's conjecture. 
This  observation   can indeed be put to use as is shown in  the examples of \cite{LNP}; the present  note   sets up   the proper theoretical  framework.

 \begin{notation}  \begin{itemize}
\item   $H^*$ denotes     Weil cohomology; 
$ \chow_*$ denotes Chow groups with $\Q$-coeffcients.
\item  The degree $d$ correspondences from $X$ to $Y$ are denoted $\Corr^d(X,Y)$.
\item  For a smooth projective manifold $X$, its Chow motive is denoted  ${\mathsf h}(X)$.
\end{itemize}
\end{notation}

\section{Motives} 
\label{sec:hlp}

Recall that a \emph{ (Chow)  correspondence of degree $k$}
from a smooth projective variety $X$ to a smooth projective variety $Y$ is a cycle class  $\chow^{\dim X+k}(X\times Y)$.
It  induces  a morphism  on Chow groups   of the same degree and on cohomology groups (of double the degree).

Correspondences can be composed and these give the morphisms in the 
 category of Chow motives. Let me elaborate briefly on this but refer to \cite{MNP} for more details.

 Precisely, \emph{an effective Chow motive} consists of a pair $(X,p)$ with $X$ a smooth projective variety and $p$
  a degree zero correspondence which is a projector, i.e., $p^2=p$.
Morphism   between motives are induced by  degree zero correspondences compatible with projectors.
 Every smooth projective variety $X$  defines a motive 
 \[
 {\mathsf h}(X)= (X,\Delta),\quad \Delta\in \chow^{\dim X}( X\times X) \text{ the class of the diagonal}
 \]
 and a morphism $f:X\to Y$ between smooth projective varieties defines a morphism $h(Y) \to {\mathsf h}(X)$ given by the transpose of the graph of $X$.  
 This procedure defines the category of effective Chow motives.
 
One can also use correspondences of arbitrary degrees provided one uses triples $(X,p,k)$ where   $p$ is  again a projector, but a morphism
$f: (X,p,k) \to (Y,q,\ell)$  is a correspondence of degree $\ell-k$ compatible with projectors. Such triples define the \emph{category of  Chow motives}.

It should be recalled that motives, like varieties have their Chow groups and cohomology groups:
\[
\aligned
\chow^m(X,p,k)&: = \im \left (  \chow^{m+k}(X) \mapright{p_*}   \chow^{m+k}(X)  \right)  ,\\
H^m(X,p,k )&:=      \im \left (  H^{m+2k}(X) \mapright{p_*}   H^{m+2k}(X)  \right) .
\endaligned
\]

Kimura \cite{Kim} has introduced the concept  finite-dimensionality for motives and he has shown that it implies the following nilpotency result.
\begin{equation}
\label{eqn:Nilp}
N  \in \Corr^0(M,M) \text{ with  trivial action on } H^*(M)  \implies N \text{ is nilpotent.}
\end{equation}

\section{The primitive motive} \label{sec:pm}

\subsection{Primitive cohomology}
Let  $M$ be a smooth projective variety. First we recall some pertinent facts concerning the Lefschetz theorems.
These are most easily described in terms of  $H_M= H^*(M)(\dim M)$, the rational cohomology  of $M$ centered in  degree $0$.
The Lefschetz operator $L\in \End (H_M)$ is a degree $2$  operator and  Hard Lefschetz states that
\begin{equation}
\label{eqn:HL}
L^j: H_M^{-j} \mapright{\sim} H_M^j,\quad j=0,\dots ,\dim M.
\end{equation}
Using this,  one can construct a linear map $\lambda: H_M \to H_M $ of degree $-2$  which is an inverse of $L$ on the subspace $LH_M\subset H_M$:
\begin{equation}
\label{eqn:Lambda}
\lambda\comp L =  L\comp \lambda  = 1 \text{ on  the image of  } L \implies   L  \comp \lambda   \text{ is a cohomological projector}.
\end{equation}
The Lefschetz decomposition gives  the corresponding  direct sum decomposition
\begin{equation}
\label{eqn:LD}
H_M= H^{\rm pr}_M  \oplus L H'_M,\quad  H^{\rm pr}_M  := \ker (\lambda),
 \end{equation}
 which is orthogonal with respect to the cup-product pairing.
Indeed, 
we have $ L  \comp \lambda (u)=0$ if $u$ is primitive and if $u= Lu'$ we have $ L  \comp \lambda (u)=L\comp  \lambda \comp L (u' )= L(u')=u$.  This shows that
\begin{equation}\label{eqn:PrimProj}
\pi^{\rm pr}:= \id -   L  \comp \lambda
\end{equation}
gives  a projector onto the primitive cohomology. 

\subsection{Construction of the "primitive" Chow projector}
We next explain under what conditions these projectors
 can be lifted to correspondences. 
 First note that $L \in \Corr^{1}(M,M)$. 
 
 \emph{Lefschetz' conjecture} $B(M)$ states that  there is
 a correspondence $\Lambda \in \Corr^{-1}(M,M)$  inducing $\lambda$.   
 More will be  needed, namely  a lift of $L^r\comp \lambda^r$ to a (Chow) projector.
 Since  $\Lambda$ and $L$  are not known to commute on the image of $p$, this motivates the following  variant  of the Lefschetz conjecture $B(M)$.
\begin{conjecture}  Property $B(M)^*$ holds if    for all $r\ge 1$  there are  correspondences 
$\Lambda_r$ and  $ \widetilde{\Lambda}_r$ in $ \Corr^{-r}(M,M)$  
such that  
\begin{itemize}
\item   $ L^{r} \comp \Lambda_r  \in \Corr^0(M,M)$ is a projector inducing $L^r\comp \lambda^r$ in cohomology.
\item  $\widetilde{\Lambda _r} \comp L^r  \in \Corr^0(M,M)$ is a projector inducing  $\lambda^r\comp L^r$ in cohomology.
\end{itemize}
\end{conjecture}
\begin{lemma}  \label{lem:BNandBM*} If $\mathsf{h}(M)$ is finite dimensional, then $B(M)$ implies  $B(M)^*$.
\end{lemma}
\proof   
I shall follow  the proof of \cite[Lemma 5.6.10]{MNP} in detail. First I shall construct $\Lambda_r$.
 Let $e= L^r \comp \Lambda^r   \in \Corr^0(M,M)$.  Since this is a cohomological projector, 
\eqref{eqn:Nilp} implies that  $e^2-e$ is nilpotent, say $(e^2-e)^N=0$. 
Introduce
 \[
\aligned
E:=  (1 - (1-e)^N)^N  &= \left(  P(e)\cdot e \right)^N ,\quad (P \text{ some polynomial})\\
&= e^N \cdot  P(e)^N  \\
&=  L^r  \comp  \Lambda^r  \cdot e^{N-1}\cdot P(e)^N .
\endaligned
\]
In cohomology this induces the same operator as $e$.
One has
\[
E=(1 - (1-e)^N)^N= 1 + \sum_{j=1}^N (-1)^j  {N \choose j}(1-e)^{jN}
\]
and so,  since $E= e^N \cdot  P(e)^N =e^N \cdot P(e)^N$,  for some polynomial $Q$ one has  
\[
\aligned
E\comp E &= E \comp (1 + \sum_{j=1}^N (-1)^j  {N \choose j}(1-e)^{jN} )\\
                 & = E +   P(e)^N \comp e^N \comp (1-e)^N \comp Q(e)\\
                 & = E  \quad \text{ (since   $e^N \comp (1-e)^N=0$)}.
                 \endaligned
                 \]   
This is thus   a projector inducing the same operator as $e$ in cohomology.  
Now  set $\Lambda_r:=  \Lambda^r  \cdot e^{N-1}\cdot P(e)^N$.  By construction   $E=L^r\comp \Lambda_r$ 
induces  the same operator as $e$ in cohomology, i.e. the operator $L^r \comp \lambda^r  $.

To show the second claim, exchange the order of $L^r$ and $\Lambda^r$.
\qed \endproof
\medskip 
For $r=1$,   this yields:
\begin{corr} Suppose that Lefschetz' conjecture $B(M)$ holds and that  the Chow motive $\mathsf{h}(M)$  is Kimura finite dimensional.
 There   is a correspondence  $\Lambda_1 \in \Corr^{-1}(M,M)$   such that $\Pi ^{\rm pr} :=\Delta_M- L \comp \Lambda_1 $ is a projector
 inducing  the projector  $\pi^{\rm pr}$  (see \eqref{eqn:PrimProj}) in cohomology.
 \end{corr}
 
 \begin{corr} For any  smooth projective  variety $M$ for which  $B(M)$ holds and with $\mathsf{h}(M)$ finite dimensional, there is a motive $(M,\Pi^{\rm pr})$ with $H^k(M,\Pi^{\rm pr})= H^k_{\rm pr}(M)$, the primitive cohomology.
 \end{corr}

\section{The variable and fixed  motive}
\subsection{Construction of the projectors}
Let  $i: X  \into  M$ be a  $d$-dimensional  smooth complete intersection of $r$ hypersurfaces.  Note that
the  graph $\Gamma_i \in X \times M$  of $i$ induces the Lefschetz correspondence   (we are ignoring multiplicative constants here)
\[
L ^r= i_*\comp i^*\in \Corr^{r}(M,M).
\]
 Set
\[
p_r:= L^r \comp \Lambda_r   \in \Corr^0(M,M),
\]
which by construction (cf.  Lemma~\ref{lem:BNandBM*})  is a projector.
One has
\begin{lemma}  \label{lem:FandVProjs}
Assume $B(M)$  and that $h(M)$ is finite dimensional.
Then the correspondences 
\[
\pi^{\rm fix}:= i^* \comp \Lambda_r \comp p_r \comp i_*  \in \Corr^0(X,X)
\]
 and 
 \[
 \pi^{\rm var}:=\Delta_X - \pi^{\rm fix}
 \]
 are   commuting projectors. 
\end{lemma}
\proof   It suffices to show that $\pi^{\rm fix}$ is a projector.
Then
\[
\aligned
(\pi^{\rm fix})^2 &= i^*  \comp \Lambda_r \comp  p_r   \comp L^r \comp \Lambda_r \comp p_r  \comp  i_*\\
                        &= i^* \comp\Lambda_r \comp  p_r ^3 i_*\\
                         &= i^*   \comp\Lambda_r \comp p_r i_*\\
                         & =\pi^{\rm fix}.  \hfill \qed
                        \endaligned
     \]
\endproof

 \subsection{Cohomological action} 
The inclusion $i:X\into M$  induces maps  $i^* :H^*(M)  \to  H^*(X) $  of degree $0$  and $i_* : H^*(X)\to H^*(M)$ of degree $2r$ with   $ i_*\comp i^*=L^r$ and $i^*\comp i_*= (L|X)^r$.

\begin{lemma}  \label{lem:vanproj}  For the action on  $H^d(X)$ one  has  $ p_r \comp i_*=i_*$ and  $\pi^{\rm fix}$ induces the projector 
$i^*\comp \lambda^r \comp  i_* $.
\end{lemma}
\proof    By definition of the fixed and variable cohomology \eqref{eqn:FixAndVarCoh}, one has   
\[
i_* H^d(X)= i_* H^d_{\rm fix}(X)= i_*\comp i^* H^{d}(M)= L^r H^{d}(M).
\]
 Since 
by equality  \eqref{eqn:Lambda},  $L$ and $\lambda$ are inverses 
on the image of $L $,   in cohomology one has  
  $p_r \comp i_* =  L^r \comp \lambda ^r  \comp i_*=i_*$ and  $\pi^{\rm fix}=  i^*\comp \lambda^r\comp i_*$. 
\qed \endproof

\begin{corr}  \label{cor:FandVinCoh}  The cohomological projectors $\pi^{\rm fix}$ and $\pi^{\rm var}$ induce projection onto the fixed and variable cohomology.
\end{corr}
\proof  Let $x\in H^d(X)$.
Then $\pi^{\rm fix}( x)=  i^* (\lambda^r\comp i_*x)  \in H^d_{\rm fix}(X)$. Since
  $i_*( x-  i^* \lambda^r i_*x) = i_*x- L^r\lambda^r i_*x= i_*x- i_*x=0$,   one has  $x-  \pi^{\rm fix} x \in H^d_{\rm var}(X)$. The result follows because of the
  direct sum decomposition   \eqref{eqn:FixAndVarCoh2}. 
\qed\endproof

 \subsection{The motives} 
 
 Now define the \emph{fixed} and  \emph{variable submotive} of $X$ by means of 
\[
\mathsf{h}(X)^{\rm fix}= (X,  \pi^{\rm fix}),\quad \mathsf{h}(X)^{\rm var}=(X, \pi^{\rm var}).
\]
Then,  Lemma~\ref{lem:FandVProjs} and Corollary~\ref{cor:FandVinCoh}  can be  summarized as follows.

\begin{prop}  \label{prop:FixAndVarMots} Let $M$ be a smooth projective manifold
for which $B(M)$ holds and suppose that  $\mathsf{h}(M)$ is finite dimensional. Let $X\subset M$ be a  smooth $d$-dimensional complete intersection.
Then $\pi^{\rm fix}$ is a projector inducing in cohomology projection onto the fixed part of the cohomology and
$
\pi^{\rm var}$ is a projector  commuting with $\pi^{\rm fix}$ and inducing projection on the variable cohomology. There is a direct sum splitting of
motives
\[
{\mathsf h}(X)= {\mathsf h}(X)^{\rm fix}\oplus {\mathsf h}(X)^{\rm var}.
\]
\end{prop}

\begin{rmk} Let $X$ be a surface. Then   \cite[\S 6.3]{MNP}  there is a self dual Chow-Lefschetz decomposition of the diagonal
\[
\Delta= \pi_0 + \pi_1 +  \underbrace{\pi_2^{\rm alg}+\pi_2^{\rm tr}}_{\pi_2}+  \pi_3+  \pi_4.
\]
This decomposition is compatible with the splitting into variable and fixed motives. This is because one has a splitting
\begin{equation} \label{eqn:Alg}
\pi_2^{\rm alg}= \pi_2^{\rm alg, fix}+ \pi_2^{\rm alg,  var}.
\end{equation}
Indeed, the construction of the projector $\pi_2^{\rm alg}$ as given in loc. cit. proceeds by first taking an orthogonal  basis for the algebraic classes  of $X$, say
$d_1,\dots, d_\rho$  with $\pi_1(d_j)=0$ for $j=1,\dots,\rho$, and then one sets
\[
\pi_2^{\rm alg} =  \sum_{i=1}^\rho  \frac{1}{ d_i^2} d_i \times d_i \in \Corr^0(X,X)
\]
Since the splitting in variable and fixed parts is an orthogonal splitting, the splitting \eqref{eqn:Alg} follows. One then puts $\pi_2^{\rm tr}= \pi_2- 
\pi_2^{\rm alg}$ and hence, defining $\pi_2^{\rm tr, var}:= \pi^{\rm var}- \pi_2^{\rm alg,var}$ and $\pi_2^{\rm tr, fix}:=\pi_2^{\rm tr}- \pi_2^{\rm tr, var}$, one gets
a refinement of the above Chow-Lefschetz decomposition
\[
\Delta= \pi_0 + \pi_1 +  \underbrace{\pi_2^{\rm alg,fix}+\pi_2^{\rm tr,fix}}_{\pi^{\rm fix}_2}+   \underbrace{\pi_2^{\rm alg,var}+\pi_2^{\rm tr,var}}_{\pi^{\rm var}_2}+  \pi_3+  \pi_4.
\]

\end{rmk}

Proposition~\ref{prop:FixAndVarMots}    asserting  the splitting into variable and fixed motives 
has  the following consequence  which states that the characterization for fixed and variable cohomology has a motivic  analog:
\begin{lemma} Same assumptions as before. \\
{\rm 1.}  For $k\le d$ we have 
\[
\chow_k (\mathsf{h}(X)^{\rm var}) = \ker (p_r \comp i_* : \chow_k(X)\to \chow_k(M)).
\]
{\rm 2.}
We have a  surjective morphism 
\[
i^*:\chow_{k+r}(M)  \onto    \chow_k (\mathsf{h}(X)^{\rm fix}) .
\]
\label{lem:AltFixAndVar}
\end{lemma}
 \proof
1. By definition the left hand side consists of cycles of the form $y=z- i^*\Lambda_r   p_r    i_* z$ for some $z\in \chow_k (X)$. Clearly,   if $p_r\comp i_*u=0$,
$u$ is of this form and conversely, if $y$ is of this form, we have  $i_*y= i_*z - i_*i^* \Lambda _r  p_r \comp i_* z= i_*z - L ^r \comp \Lambda _r  p_r \comp i_* z=
i_*z -   p_r \comp i_* z$ since $p_r$ is a projector
and applying $p_r$ this vanishes.  
\\
2. This follows since the "fixed" cycles are all in the image of $i^*$.
 \qed\endproof


\subsection{Variants with group actions} \label{sec:VarGrps}
 There are variants with group actions as follows. Suppose that a finite group $G$ acts on $M$ and that $X$ is invariant under the action of $G$. In particular, $g$
commutes with $i$ and with $L_X$ and $L_M$.
Let $\Gamma_g$ be the graph of the action of $g$ on  $X$. For $\chi = \sum_g \chi(g)\cdot g \in \Q[G]$ we set
\[
\pi_\chi := \frac{1}{|G|} \sum \chi(g) \Gamma_{g}.
\]
This is a projector and defines the motive $(X,\pi_\chi)$. Moreover,
$\pi_\chi \comp \pi^{\rm fix}=  \pi^{\rm fix}\comp  \pi_\chi$ and hence, $ \pi^{\rm fix}$ also commutes with $ \pi^{\rm var}$ and both
$\pi^{\rm fix}  \comp \pi_\chi$ and $\pi^{\rm fix}  \comp \pi_\chi$ are projectors.
For any $\Q$-vector space on which $G$ acts, setting
\[
V^\chi:= \sett{x\in V}{ {g(x)= \chi(g) x \text{ for al } g\in G}},
\]
one has
\[
H^k (X,\pi_\chi) = H^k(X)^\chi.
\]
Since $X\subset M$ is left invariant by the $G$-action, the variable and fixed motives are $G$-stable and one sets
\[
\mathsf{h}(X,\pi_\chi)^{\rm fix}:=   (X,  \pi^{\rm fix}  \comp \pi_\chi),\quad   \mathsf{h}(X,\pi_\chi)^{\rm var }:=   (X,  \pi^{\rm var}  \comp \pi_\chi) .
\]

\medskip

  \begin{nonumberingt} Thanks to Robert Laterveer and Jaap Murre for their  interest and   suggestions. 
\end{nonumberingt}

\end{document}